\documentclass[11pt]{article}
\usepackage[utf8]{inputenc}
\usepackage{amsthm}
\usepackage{graphicx}
\usepackage{authblk}
\usepackage[T1]{fontenc}

\newtheorem{thm}{Theorem}[section]

\newtheorem{cor}[thm]{Corollary}

\newtheorem{lem}[thm]{Lemma}

\input epsf
\date{}

\begin{document}

\title{\textbf{Some New Results in the Alcuin Number of Graphs}}
\author[1]{Abbas Seify \thanks{abbas.seify@srttu.edu}}
\author[2]{Hossein Shahmohamad \thanks{hxssma@rit.edu}}

\affil[1]{\small{Department of Mathematical Sciences, Shahid Beheshti University, Evin, Tehran, Iran}}
\affil[2]{\small{School of Mathematical Sciences, R.I.T, Rochester, NY 14623}}

\maketitle

\begin{abstract}
We prove some results concerning Alcuin number of graphs. First, we classify graphs which have unique minimum vertex cover. Then we present two necessary conditions for a graph to be of class two and show why one of them (condition on common neighbors) is sufficient as well. By using this classification theorem, we prove some results about class one and class two graphs such as cartesian product of graphs and classification of regular graphs.
\end{abstract}

\section{Introduction}
Let $v\in V (G)$. We denote the set of all neighbors of $v$ by $N(v)$ and for $X \subseteq V (G)$ we define $N(X) = \cup _{x \in X}N(x)$. Also, we denote the neighbors of $X$ in $S$ by $N_S(X) = N(X) \cap S$. An {\it independent set} is a set of vertices in a graph, no two of which are adjacent. The {\it independence number}, $\alpha(G)$, is the size of greatest independent set in $G$. A {\it vertex cover} is a set of vertices such that every edge in graph has an end in it. The {\it vertex cover number} is the size of smallest vertex cover in $G$ and is denoted by $\beta(G)$.The {\it girth} of a graph $G$, is denoted by $g(G)$, is the length of a shortest cycle contained in the graph. If the graph does not contain any cycles, its girth is defined to be infinity.\\
Let $G$ and $H$ be two simple graphs. The {\it cartesian product} of $G$ and $H$ is denoted by $G\times H$ and is a graph such that $V (G\times H) = \{(u,v)| u\in G , v \in H\}$ and two vertices $(u,v)$ and $(x,y)$ are adjacent if and only if $u = x$ and $v$ is adjacent to $y$ or we have $v = y$ and $u$ is adjacent to $x$. For $n = 0$ we define $Q_0 = K_1$ and for $n \geq 1$, we define {\it $n$-dimensional hypercube} as $Q_n = Q_{n-1} \times K_2$.\\\\
In [2], Alcuin’s river crossing problem and its history is fully explained:
\begin{quote}
A man had to transport to the far side of a river, a wolf, a goat, and a bundle of cabbages. The only boat he could find was one which would carry only himself and one of them. For that reason he sought a plan which would enable them all to get to the far side unhurt.
\end{quote}
We consider a graph $G = (V,E)$ where for each object we have a vertex in $V(G)$ and two vertices are adjacent if and only if corresponding objects are conflicting. The aim is to transport objects to the far side of river unhurt i.e. transporting the vertices of graph to the far side such that in each step of transportation, there exists no edges between the vertices which are in the same side of river. We assume that in the starting step objects are in the left side. For example, in Alcuin river crossing puzzle if we consider $w$ as wolf, $g$ as goat and $c$ as cabbage, then figure 1.1 shows a feasible schedule for problem.\\
\begin{figure}
\centering{\includegraphics[width=10cm]{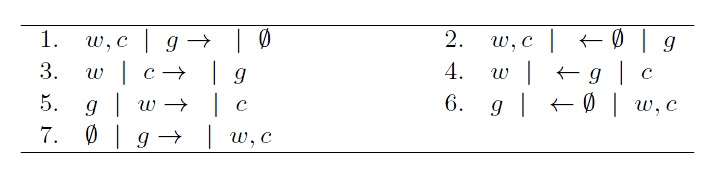}} 
\caption{Feasible Schedule for Alcuin Puzzle}
\end{figure}
The {\it Alcuin Number} of a graph $G$ is denoted by $c(G)$ and is defined as the least capacity of a boat that can successfully handle this problem. In [1],[2],[3] the following has been proved.
\begin{thm}
For any graph $G$ we have $\beta(G) \leq c(G) \leq \beta(G) + 1$.
\end{thm}
\begin{proof}
In the first step we should have at least $\beta(G)$ of vertices in boat. This shows that $c(G) \geq \beta(G)$. Consider a boat of capacity $\beta(G)+1$. Put a minimum vertex cover $C \subseteq V (G)$ in boat. Now, boat has one empty capacity. Put the vertices of $V-C$ in boat and transport them to the right side of river one by one. By repeating this prosedure we can transport $V-C$ to the right side. In final step, drop off $C$ in the right side. This is a feasible schedule for $G$ by a boat of capacity $\beta(G) + 1$ and completes the proof.
\end{proof}
According to this theorem, graphs are divided into two classes. We say a graph $G$ is of calss one if $c(G) = \beta(G)$ and of class two if $c(G) = \beta(G)+1$. Our goal is to determine class one and class two graphs. In section 2 we prove the following theorem.
\begin{thm}
\textbf{(Classification Theorem)} Let $G$ be a graph and $C\subseteq V (G)$ be a minimum vertex cover of $G$. Then $G$ is of class two if and only if for any two independent sets $S,T \subseteq C$ we have $|N_{V-C}(S) \cap
N_{V-C}(T)| > |S| + |T|$.
\end{thm}
Alcuin's river crossing problem differs significantly from other medieval puzzles, since it is neither geometrical nor arithmetical, but purely combinatorial [2]. In [3], the Ferry Problems, which may be viewed as generalizations of the classical wolf-goat-cabbage puzzle, are considered. The Ferry Cover Problem, where the objective is to determine the minimum required boat size to safely transport $n$ items represented by a graph $G$, is studied in [3]. A close connection with vertex cover which leads to hardness and approximation results is also studied in [3].

\section{Classification Theorem}
In this section we classify graphs which have unique minimum vertex cover and find some necessary conditions on graphs to be of class two. Finally, we prove the classification theorem, which presents a necessary and sufficient condition for graphs to be of class two. In [2], the following has been proved.

\begin{thm}
Let $G=(V,E)$ be a graph which has two distinct independent sets $S_1,S_2 \subseteq V(G)$ of maximum size (or, equivalently, two distinct vertex covers of minimum size). Then $G$ is of class one.
\end{thm}
In the following theorem, we classify the graphs which have unique minimum vertex cover. So, we can find a necessary condition on graphs to be of class two.

\begin{thm}
Let $G =(V,E)$ be a graph and $C \subseteq V (G)$ be a minimum vertex cover of $G$. Then $C$ is unique minimum vertex cover in $G$ if and only if for every independent set $A\subseteq C$ we have $|N_{V-C}(A)| > |A|$.
\end{thm}
\begin{proof}
First, suppose that $C$ is unique minimum vertex cover of $G$. Note that for every independent set $A\subseteq C$  we have $|N_{V-C}(A)|\geq |A|$. Because, if there exists an independent set $A\subseteq C$ such that $|N_{V-C}(A)| <|A|$, then by replacing $A$ with $N_{V-C}(A)$ we obtain a vertex cover $C'\subseteq V (G)$ such that $|C'|< |C|$, which is a contradiction.\\
So, it sufficies to show that $|N_{V-C}(A)| \neq |A|$. If not, by the similar way, we can find another minimum vertex cover of $G$ of cardinality $|C|$ which is a contradiction, since $C$ is unique minimum vertex cover of $G$.\\
For converse, suppose that the condition holds. By contrary, let $G$ have two distinct minimum vertex covers $C$ and $D$.\\
Suppose that $C \cap D = \emptyset$. Then $C$ and $D$ are independent sets and $N_{V-C}(C) \subseteq D$. Because, if there exists an edge $e = uv$ such that $u\in C$ and $v\in V-(C \cup D)$, then $e$ is not covered by $D$ which is a contradiction. Now, let $A = C$. Then we have $|N_{V-C}(A)|\leq |A|$ which is a contradiction.\\
So, we may assume that $C \cap D \neq \emptyset$. Then $C-D$ and $D-C$ are independent sets and by the similar way to the first case, we can show that $N_{V-C}(C-D) \subseteq D-C$. Therefore, we have $|N_{V-C}(C-D)| \leq |D-C| = |C-D|$, which is a contradiction and this completes the proof.
\end{proof}
By using this theorem, we can find a necessary condition for class two graphs.

\begin{cor}
Let $G$ be of class two and $C\subseteq V(G)$ be a minimum vertex cover of $G$. Then for every independent set $A\subseteq C$ we have $|N_{V-C}(A)|>|A|$.
\end{cor}
This condition is not sufficient. For example, consider $K_{1,2}$. Then it satisfies the condition, but we have $c(G) = \beta(G) = 1$, which shows the graph is of class one. In following, we present a better condition.

\begin{thm}\label{independent}
Let $G$ be of class two and $C \subseteq V (G)$ be a minimum vertex cover of $G$. Then
for every independent set $A \subseteq C$ we have $|N_{V-C}(A)|>2|A|$.
\end{thm}

\begin{proof}
By contrary, suppose that there exists an independent set $A\subseteq C$ such that $|N_{V-C}(A)|\leq 2|A|$. Then, following algorithm shows that $c(G) = \beta(G)$.
\begin{enumerate}
\item Put $C$ in boat and go to the right side. Leave $A$ in the right side and return to the left side.
\item By $|A|$ empty places in boat, in each return to the left side, transfer $|A|$ elements of the set $(V-C)-N_{V-C}(A)$ to the right side and drop off them.
\item Now, first transfer $|A|$ elements of $N_{V-C}(A)$ to the right side and drop off them. Then put $A$ in boat, return it to the left side and leave it in left side and transfer remaining elements of $N_{V-C}(A)$ to the right side of river (Note that this is possible because $|N_{V-C}(A)| \leq 2|A|$).
\item Finally, return to the left side, put $A$ in boat, go to the left side and drop off $C$ in the right side.
\end{enumerate}
This algorithm transports $G$ by a boat of capacity $\beta(G)$. So, $G$ is of class one and this completes the proof.
\end{proof}
This theorem has an interesting and simple corollary about claw-free graphs.

\begin{cor}
Let $G$ be a claw-free graph. Then $G$ is of class one.
\end{cor}
\begin{proof}
Let $G$ be of class two and $C\subseteq V(G)$ be a minimum vertex cover of $G$. Consider $v\in C$. Then \ref{independent} implies that $v$ has at least three distinct neighbors in $V-C$ such as $x,y,z$. Since $V-C$ is an independent set, $v$ and $x,y,z$ induce a claw and this completes the proof.
\end{proof}
There are natural questions: Is this condition sufficient? Is there a natural number $k$ such that if $|N_{V-C}(A)|>k|A|$ for all independent sets $A\subseteq C$, then graph is of class two? Both questions answers are negetive. Let $G = (V,E)$ such that $V =\{v_1, v_2, u_1, u_2, u_3, \ldots , u_{2k+1}\}$ and $E = \{v_1u_1, \ldots , v_1u_{k+1}, v_2u_{k+1}, \ldots , v_2u_{2k+1}\}$. Then $C = \{v_1, v_2\}$ is a minimum vertex cover of $G$. Also, it can be checked easily that $G$ is of class one and satisfies the conditions. Now, we present the classification theorem.

\begin{thm}\label{classification}
\textbf{(Classifcation Theorem)} Let $G$ be a graph and $C \subseteq V (G)$ be a minimum vertex cover of $G$. Then $G$ is of class two if and only if for any two independent sets $S,T \subseteq C$ we have $|N_{V-C}(S) \cap N_{V-C}(T)| > |S| + |T|$.
\end{thm}
We will use the structure theorem in [2] to prove the classification theorem. The structure theorem is as follows.
\begin{thm}
\textbf{(Structure Theorem)} A graph $G = (V;E)$ possesses a feasible schedule for a boat of capacity $b \geq 1$ if and only if there exist five subsets $X_1, X_2, X_3, Y_1, Y_2$ of $V$ that satisfy the following four conditions:
\begin{enumerate}
\item The three sets $X_1,X_2,X_3$ are pairwise disjoint.Their union $X = X_1 \cup X_2 \cup X_3$ forms an independent set in $G$.
\item The (not necessarily disjoint) sets $Y_1,Y_2$ are nonempty subsets of the set $Y=V-X$, which satisfies $|Y| \leq b$.
\item $X_1 \cup Y_1$ and $X_2 \cup Y_2$ are independent sets in $G$.
\item $|Y_1| + |Y_2| \geq |X_3|$
\end{enumerate}
\end{thm}
Note that if $G$ is of class one, then $Y$ is a minimum vertex cover and $X$ is a maximum independent set in $G$. Now, by using the structure theorem we can prove classification theorem.

\begin{proof}
First, suppose that there exist independent sets $S,T \subseteq C$ such that $|N_{V-C}(S) \cap N_{V-C}(T)| \leq |S|+|T|$. Then the following transportation algorithm shows that $G$ is of class one.
\begin{enumerate}
\item Transport $C$ to the right side of river, drop off $S$ in right side and return to the left side.  
\item  By $|S|$ empty places in boat, transport $|S|$ elements of the set $(V-C) -N_{V-C}(S)$ to the right side. By repeating this procedure, transport the set $(V-C)-N_{V-C}(S)$ to the right side and return to the left side.
\item  Now, by using $|S|$ empty places in boat, put $|S|$ elements of the set $N_{V-C}(S)\cap N_{V-C}(T)$ in boat, transport them to the right side and drop off them. Then put $S$ in boat and return to the left side.
\item  In left side, drop off $T$ and put remaining elments of $N_{V-C}(S) \cap N_{V-C}(T)$ in boat and transport them to the right side and leave them. (Note that this is possible, because we have $|N_{V-C}(S) \cap N_{V-C}(T)| \leq |S| + |T| $)
\item  Now, by $|T|$ empty places in boat, transfer the remaining items to the right side. In final movement, put $T$ in boat and drop off $C$ in right side.
\end{enumerate}
For converse, by contrary suppose that $G$ is of class one. Then in structure theorem $Y$ is a minimum vertex cover, $X$ is a maximum independent set and there exists a feasible schedule for a boat of capacity $b = \beta(G)$. By condition three in the structure theorem we have $X_1 \cup Y_1$ and $X_2 \cup Y_2$ are independent sets. So, $Y_1$ and $Y_2$ are independent sets in vertex cover and therefore common neighbors of $Y_1$ and $Y_2$ in $V-Y$ are in $X_3$. So, we have $|X_3| \geq |N_{V-Y}(Y_1) \cap N_{V-Y}(Y_2)| > |Y_1| + |Y_2|$, which is a contradiction, since condition four in the structure theorem indicates that $|Y_1| + |Y_2| \geq |X_3|$. Hence, $G$ is of class two.
\end{proof}

\section{Results}
In this section, we prove some results concerning class one and class two graphs such as girth of class two graphs, classification of regular graphs and cartesian product of graphs.

\begin{cor}
Let $G$ be of class two and $C \subseteq V (G)$ be a minimum vertex cover of $G$. Then for every independent set $S \subseteq C$ we have $|N_{V-C}(S)| > 2|S|$.
\end{cor}
\begin{proof}
In theorem \ref{classification} let $S = T$.
\end{proof}
A simple and exciting result of classification theorem is the case in which $S$ and $T$ are single vertex. This yields to the following corollary which has been proved in [2] by considering suitable subsets of $V(G)$ in the structure theorem.

\begin{cor}\label{vertex}
Let $G=(V,E)$ be a graph and there exists two (not necessarily distinct) vertices $u$ and $v$ such that they have at most two common neighbors in $V-C$. Then $G$ is of class one.
\end{cor}
Another easy corollary arises about the girth of class two graphs.

\begin{cor}\label{girth}
Let $G$ be of class two with $\beta(G) \geq 2$. Then $g(G) \leq 4$.
\end{cor}
\begin{proof}
Let $C$ be a minimum vertex cover of $G$. Consider two distinct vertices $u$ and $v$ in $C$. These vertices have at least three common neighbors in $V-C$ such as $\{x,y,z\}$. Now, consider a cycle $u,x,v,y,u$. This implies that $g(G) \leq 4$.
\end{proof}
We use this corollary to prove the following corollary which also has been proved in [1].

\begin{cor}
Let $T$ be a tree. Then $T$ is of class two if and only if $T = K_{1,n}$ for $n \geq 3$.
\end{cor}
\begin{proof}
Let $T$ be a tree. Note that if $\beta(G) \geq 2$, then \ref{girth} implies T is of class one. So, we may assume that $\beta(G) = 1$. Therefore $T = K_{1,n}$ and the classification theorem completes the proof.
\end{proof}
We can find a better result about girth of regular graphs.
\begin{thm}
Let $G$ be a regular graph of class two. Then $g(G)=3$.
\end{thm}
\begin{proof}
Let $G$ be an $r$-regular graph of class two. Then \ref{vertex} implies that $r \geq 3$ and this yields that $\beta(G) \geq 2$. Suppose that $G$ is not bipartite. Then there exist two vertices $u,v$ in minimum vertex cover of $G$ such that $uv \in E(G)$. Since $G$ is of class two, these two vertices have at least three common neighbors in $V-C$ and therefore $G$ has a triangle. So, $g(G) = 3$.\\
So, we may assume that $G=(A,B)$ is bipartite. Since $G$ is regular we have $|A|=|B|$. On the other hand, bipartite regular graphs have perfect matching. So, we have $\beta(G) \geq |A|$. Also, $A$ is a vertex cover of $G$. This implies that $\beta(G) \leq |A|$. Therefore $\beta(G) = |A|$ and $A$ (or $B$) is a minimum vertex cover of $G$. Now, first transfer $A$ to the right side and then return and transfer $B$. So, $c(G) = \beta(G)$ which implies that $G$ is of class one which is a contradiction and this completes the proof.
\end{proof}
Next, we prove a theorem that classifies the cartesian product of graphs.

\begin{thm}
Let $G$ and $H$ be two graphs such that both of them are not trivial. Then $G \times H$ is of class two if and only if one of them is of class two and another is $K_1$.
\end{thm}
\begin{proof}
Let $G$ and $H$ satify the conditions. Then $G \times H$ is of class two. For converse, by contrary, suppose that conditions don't hold. If $G$ or $H$ is $K_1$ and the other is of class one, then $G \times H$ is a of class one which is a contradiction.\\
Now, suppose that both of $G$ and $H$ have at least two vertices. If one of them is an empty graph, then $G \times H$ is disconnected and has at least two nontrivial components. So, there exists one vertex of minimum vertex cover in each component. These vertices don't have any common neighbors. Therefore \ref{vertex} implies that $G$ is of class one and this completes the proof.\\
So, we may assume that $G$ and $H$ are nontrivial and have at least two vertices. First, suppose that $max\{\beta(G) , \beta(H)\} \geq 2$. We may assume that $\beta(G) \geq 2$. Since, $|V (H)| \geq 2$ we can find two vertices $(u,v)$ and $(x,y)$ in minimum vertex cover of $G \times H$ such that $u \neq x , v \neq y$ (to find such vertices, it suffices to consider two copies of $G$ in $G \times H$). Now, consider these vertices. Clearly, they have at most two common neighbors and therefore $G \times H$ is of class one which is a contradiction.\\
So, we may assume that $max\{\beta(G),\beta(H)\} = 1$. Suppose that $V (G) = \{v_1, \ldots ,v_n\}$ and $V (H) = \{u_1,\ldots,u_m\}$ and mininmum vertex covers of $G$ and $H$ are $\{v_1\}$ and $\{u_1\}$ respectively. Also, we may assume that $m , n \geq 3$, since otherwise $G \times H$ is $C_4$ which is of class one. So, minimum vertex cover of $G \times H$ is $\{(v_1,u_1),(v_1,u_2),\ldots,(v_1,u_m),(v_2,u_1),(v_3,u_1),\ldots,(v_n,u_1)\}$. Consider $(v_1, u_2),(v_2,u_1)$. These two vertices have at most two common neighbors and therefore $G \times H$ is of class one which is a contradiction and this completes the proof.
\end{proof}
Now, we can conclude the following.
\begin{cor}
The $n$-dimentional hypercube $Q_n$ is of class one for all $n \geq 1$.
\end{cor}
Now, we prove a result concerning class of regular graphs. We will use the following lemma in the next theorem.

\begin{lem}\label{lem}
Let $G$ be an $r$-regular graph ($r \geq 1$) of order $n$ and $C \subseteq V (G)$ be a mninmum vertex cover of $G$. Then $|C| \geq max\{\frac{n}{2} , r\}$.
\end{lem}
\begin{proof}
Let $G$ be an $r$-regular graph. Then $\alpha(G) \leq \frac{n}{2}$ and also we have $\alpha(G) + \beta(G) = n$ which implies that $\beta(G) \geq \frac{n}{2}$ . Also, note that since $r \geq 1$, we can consider $v \in V-C$. Then all neighbors of $v$ are in $C$ which implies that $|C| \geq r$ and this completes the proof.
\end{proof}

\begin{thm}
Let $G$ be an $r$-regular graph and $r \in \{2,3,4,5\}$. Then $G$ is of class one.
\end{thm}
\begin{proof}
Case $r = 2$ is clear. Let G be a cubic graph and $C \subseteq V (G)$ be a minimum vertex cover of $G$. Then \ref{lem} implies that $|C| \geq 2$. If $C$ is an independent set, then $G$ is bipartite and therefore of class one. So, suppose that there exists an edge $e = uv$ such that $u,v \in C$. Since $G$ is cubic, $u$ and $v$ have at most two neighbors in $V-C$ and \ref{vertex} implies that $G$ is of class one.\\\\
Let $G$ be a 4-regular graph. If $C$ is an independent, then $G$ is bipartite and we are done. So, consider $u,v \in C$ such that $uv \in E(G)$. If $u$ and $v$ have at most two common neighbors in $V-C$, then we are done. So, suppose that $u$ and $v$ have three common neighbors in $V-C$. Since $G$ is 4-regular, this implies that $u$ has no neighbor in $C$ except $v$ and also similar is true for $v$. Now \ref{lem} implies that $|C| \geq 3$ and there exists $z \in C -\{u,v\}$. Let $S = \{u,z\}$ and $T = \{v\}$. Clearly, $S$ and $T$ are independent sets. If $G$ is of class two then classification theorem yields that $|N_{V-C}(S) \cap N_{V-C}(T)| \geq 4$ which is impossible, since $v$ has three neighbors in $V-C$. Therefore $G$ is of class one.\\\\
Let $G$ be a 5-regular graph. Then \ref{lem} implies that $|C| \geq 5$. Similar to the previous case we may assume that there exists vertices $u$ and $v$ in $C$ such that $uv \in E(G)$. If $G$ is of class two, then $u$ and $v$ have at least three common neighbors in $V-C$.\\
If there exists $z \in C- \{u,v\}$ which is adjacent to $u$, then $deg_{V-C}(u) = 3$ and we can find some vertex $t$ such that $t$ is not adjacent to $v$ (This is possible since $|C| \geq 5$). Now, let $S = \{v,t\}$ and $T = \{u\}$. Then, if $G$ is of class two, classification theorem implies that $|N_{V-C}(S) \cap N_{V-C}(T)| \geq 4$ which is impossible, since $deg_{V-C}(u) = 3$.\\
So, we may assume that $u$ and $v$ have no other neighbor in $C$. Consider $z \in C-\{u,v\}$. Since induced subgraph on $C-\{z\}$ is not complete, so there exists an independent set $S \subseteq C-\{z\}$ such that $|S| \geq 2$. Let $T = \{z\}$. Then classification theorem implies that $deg_{V-C}(z) \geq 4$. This implies that for every $c \in C$ we have $deg_C(c) \leq 1$. So, we can find an independent set $S \subseteq C$ such that $|S| \geq 3$. It is clear that both of $u$ and $v$ are not in $S$. Suppose that $u$ is not in $S$. Let $T = \{u\}$. Then classification theorem implies that $|N_{V-C}(S) \cap N_{V-C}(T)| \geq 5$, which is impossible since $deg_{V-C}(u) = 4$.
\end{proof}
There are some open questions about the problem as follows.\\
\textbf{Question 1:} Are all regular graphs of class one?\\
\textbf{Question 2:} Let $G$ be a graph of order $n$. Does the probability of being class two tend to zero as $n \rightarrow \infty$? i.e. Is it true that almost all graphs are of class one?

\section*{Acknowledgements}
A recent work on the above problems is [1] which is a beautifully written play. It presents the Alcuin's river crossing problem and the ferry problems for the readers and then generalizes the results and offers numerous new problems and some conjectures. The authors of this paper would like to take the time to express their appreciation for the work done by Mehdi Behzad in [1] for exposing more mathematics to the general public. This play script, first written in Farsi, was coauthored by M. behzad and N. Samini and was published in 2011. Then its English translation by M. Behzad was edited skillfully by Professor Cheryl E. Praeger[1].\\
Here is a part of the introduction written by Professor Praeger, ``... If the mathematical developement sketched in the play had really happened in the decades and centuries following publication of the wolf, goat and cabbage riddle, then the important mathematical disciplines of Discrete Mathematics (Combinatorial and Graph Theory), Optimization (Linear and Integer Programming) and Operation Research may have appeared a millennium
earlier...''\\
Professor Dr. Beutelspacher indicates that ``... The project introduced by Professor Mehdi Behzad is an outstanding example of modern way of popularizing mathematics ...''


\begin{thebibliography}{9}

\bibitem{1} Mehdi Behzad, Naghme Samini, The Legend of The King and The Mathematician, Candle and Fog Publishing, 2013
\bibitem{2} Peter Csorba, Cor A. J. Hurkens, Gerhard J. Woeginger, The Alcuin number of a graph and its connections to the vertex cover number , Siam J. Discrete Math, Vol. 24, No. 3, pp. 757-769, 2010
\bibitem{3} Michael Lampis, Valia Mitsou, The Ferry Cover Problem, Theory Comput Syst (2009) 44: 215-229

\end{thebibliography}
\end{document}